\documentclass[12pt]{article}
\usepackage{amsmath}
\usepackage{amssymb}
\usepackage{theorem}
\usepackage{amscd}

\newtheorem{theorem}{Theorem}

\newtheorem{prop}{Proposition}

\newtheorem{lemma}{Lemma}

\newcommand{\codim}{\mathrm{codim}\,}
\newcommand{\ddim}{\mathrm{dim}\,}
\newcommand{\cchar}{\mathrm{char}\,}

\DeclareMathOperator{\Hom}{Hom}

\newcommand{\CC}{\mathbb C}
\newcommand{\kk}{{\bf k}}

\newcommand{\Oo}{\mathcal O}

\begin{document}

\begin{center}
{\Large\bf Multiplicative structure on the Hochschild cohomology 
of crossed product algebras\\[0.7cm]}

Rina Anno \\[0.7cm]
\end{center}

\begin{abstract}
Consider a smooth affine algebraic variety $X$ over an algebraically closed field $\kk$,
and let a finite group $G$ act on $X$.
We assume that $\cchar \kk$ is greater than $\ddim X$ and $|G|$.
An explicit formula for multiplication
on the Hochschild cohomology of a crossed product $HH^*(\kk[G]\ltimes\kk[X])$
is given in terms of multivector fields on $X$ and $g$-invariant
subvarieties of $X$ for $g\in G$.
\end{abstract}

\section{Introduction}

Let $X$ be a smooth algebraic variety.
It is well-known
(\cite{Caldararu}, \cite{HKR}, \cite{Kontsevich}) that 
the groups $Ext^*_{X\times X} (\Oo_{\Delta}, \Oo_{\Delta})$ 
where $\Delta$ denotes the diagonal (these groups will be
further referred to as Hochschild cohomology $HH^*(X)$) may be interpreted in terms
of multivector fields on $X$:
$$
HH^i(X) = \bigoplus\limits_{p+q=i} H^p (X, \Lambda^q TX)
$$
as vector spaces (\cite{Kontsevich}, Thm. 8.4, \cite{Caldararu}, Corr. 4.2).
For affine $X$ this becomes an algebra isomorphism
$$
HH^*(X) \cong \bigoplus\limits_i \Gamma (X, \Lambda^i TX)
$$
where the multiplication is given by the cup product
in the left hand side, and the wedge product
in the right hand side.

Our goal is to extend this result to the situation of a smooth affine algebraic $X$
with an action of a finite group $G$. 

The author would like to thank Pavel Etingof, Vasiliy Dolgushev, Boris Shoikhet and Sergei Fironov
for useful discussions. The author would also like to express her special gratitude
to Xiang Tang for pointing out an essential gap (which led to a significant
shortening of the paper).

\section{Main results}

Let $X$ be a smooth affine algebraic variety over an
algebraically closed field $\kk$ of characteristic greater than $\ddim X$ and $|G|$.
Denote by $A= \Gamma(X, \Oo_X)$ the algebra
of regular functions on $X$, let $B = \kk[G] \ltimes A$ be crossed product algebra with multiplication
defined by $ga=a^gg$, where $a^g$ denotes the result of $g$-action on $a$.
Denote by $\pi^g$ the operator of symmetrization by $g$: if $g$ has order $k$, then
$\pi^g = \frac{1}{k}\sum\limits_{i=1}^k g^i$. Let it act on the multivectors of degree $l$
as $\pi^g\wedge\ldots\wedge \pi^g$ ($l$ components).
Let $X^g_m$ for $g \in G$ denote the $m$-th connected component of
the subvariety of $g$-invariants in $X$, let $d_{g, m}=\codim X^g_m$.

\begin{theorem}\label{isom} {\rm (cf. \cite{NPPT}, Theorem 3.11)}
$$
HH^i(B) \cong 
( \bigoplus\limits_{g \in G, m} \Gamma(X^g_m, (\Lambda^{i - \codim X^g_m} TX^g_m \otimes \Lambda^{top}(N_XX^g_m)))^G
$$
as vector spaces.
\end{theorem}

A similar result was proven in \cite{GK}, \cite{Kaledin}:
\begin{theorem}\label{sympl} {\rm (\cite{GK}, \cite{Kaledin})} 
For a complex symplectic vector space $V$ with a symplectic linear
action of a finite group $G$ one has
$$
HH^i(\CC[G]\ltimes \CC[V]) \cong 
(\bigoplus\limits_{g \in G} \Omega^{i-\codim V^g} (V^g))^G
$$
as graded vector spaces.
\end{theorem}
Here $\Omega^k(Y)$ denotes the space of differential $k$-forms on $Y$.
One can easily see that in presence of a symplectic form
there is an isomorphism
$\Omega^i(X^g_m) \cong \Lambda^{i - codim X^g_m} TX^g_m \otimes \Lambda^{top}(N_XX^g_m)$,
so Theorem \ref{sympl} is a special case of Theorem \ref{isom}.

\begin{theorem}\label{mult}
The multiplication on $HH^i(B)$ is given by
\begin{multline*}
(\sum\limits_{g\in G, m} \alpha_{g, m} \otimes \beta_{g, m})\cdot (\sum\limits_{h\in G, n} \gamma_{h, n} \otimes \delta_{h, n}) = \\
= (-1)^{d_{g, m}(j-d_{h, n})}
\sum\limits_{u\in G, k}\sum _{\substack{gh=u,\\ X^g_m\cap X^h_n\supset X^u_k}}
(\pi^u\alpha_{g, m}|_{X^u_k}\wedge\pi^u\gamma_{h, n}|_{X^u_k}) \otimes (\beta_{g, m}|_{X^u_k}\wedge\delta_{h, n}|_{X^u_k})
\end{multline*}
where $\alpha_{g, m}, \gamma_{g, m} \in \Gamma(X^g_m, \Lambda^{i - d_{g, m}} TX^g_m)$ and
$\beta_{g, m}, \delta_{g, m} \in \Gamma(X^g_m, \Lambda^{d_{g, m}}N_XX^g_m)$.
The sum in the RHS is taken over all $g$ and $h$ such that $X^{gh}_k$ is a component of $X^g_m\cap X^h_n$
for some $k, m, n$.
\end{theorem}

Another description of the multiplication in $HH^*(B)$ was given in \cite{Witherspoon}.

When $X$ is a symplectic variety, the formula simplifies in the following way:

\begin{theorem}\label{main_sympl}
Let $X$ be a symplectic variety with a symplectic form $\omega$, and let the action of
$G$ be symplectic. Then
$$
HH^i(B) \cong 
( \bigoplus\limits_{g \in G} \Gamma(X^g_m, \Lambda^{i - d_{g, m}} TX^g_m))^G
$$
as vector spaces. The multiplication is given by 
\begin{equation*}
(\sum\limits_{g\in G, m} \alpha_{g, m}) \cdot (\sum\limits_{h\in G, n} \gamma_{h, n}) =
(-1)^{d_{g, m}(j-d_{h, n})}
\sum\limits_{u\in G, k}\sum _{\substack{gh=u,\\ X^g_m\cap X^h_n \supset X^u_k}}
(\pi^u\alpha_{g, m}|_{X^u_k}\wedge\pi^u\gamma_{h, n}|_{X^u_k}).
\end{equation*}
\end{theorem}

\section{Preliminaries}

Define an $A$-bimodule $Ag$ as a submodule $Ag=\{ag\ |\ a\in A\}$ in $B$.
The following proposition was proved in \cite{DE}:
\begin{prop}\label{A_Ag} {\rm (\cite{DE}, Prop. 3)}
$$HH^i(B) \cong (\bigoplus\limits_{g \in G} H^i(A, Ag))^G$$
as vector spaces.
\end{prop}

In fact, the complexes $C^*(B,B)$ and 
$(\bigoplus\limits_{g \in G} C^*(A, Ag))^G \cong (C^*(A,B))^G$
are quasiisomorphic. The multiplication is defined on the Hochschild cochains,
hence on $(C^*(A,B))^G$. To write it down in a convenient basis one
should extend it to $C^*(A,B)$ (possibly losing cohomological properties
beyond $G$-invariants). We have a following lemma-definition:

\begin{lemma}\label{mu}
The multiplication on the Hochschild cochains $C^*(A,B)$ given by the map
$$\mu: C^i(A,B) \otimes C^j(A,B) \to C^{i+j}(A,B)$$
$$
\mu (\phi g \otimes \psi h) (a_1\otimes\ldots\otimes a_{i+j}) =
\phi(a_1\otimes\ldots\otimes a_i)\cdot\psi^g(a_{i+1}\otimes\ldots\otimes a_{i+j}) gh,
$$
where $\phi \in \Hom(A^{\otimes i}, A)$, $\psi \in \Hom(A^{\otimes j}, A)$,
on the $G$-invariant cohomology $(H^*(A,B))^G$ coincides with the multiplication
that comes from the isomorphism with $HH^*(B)$.
\end{lemma}

The proof is a direct computation.

\section{The local case}

In this section let $X$ be a linear space $V\simeq \kk^n$, so 
$\kk[V]=\kk[x_1,\ldots x_n]$ is the polynomial algebra. Let the action of the group $G$ on $X$
be linear. For $g \in G$ denote by $V^g$ the space of $g$-invariant vectors in $V$, 
and by $(V^g)^\vee$ the subspace of $V$ generated by eigenvectors of $g$
with eigenvalues different from one. The matrix of $g$ is diagonalizable,
so $V = V^g \oplus (V^g)^\vee$. 
Note that $\kk[V^g]g \simeq \kk[V^g]$ as an $A$-bimodule (compare with Lemma \ref{simeq}).

\begin{prop}\label{loc_isom}
\begin{enumerate}
\item $H^i(\kk[V],\kk[V]g) \cong \Lambda^{i - d_g} V^g \otimes  \Lambda^{d_g}(V^g)^\vee \otimes \kk[V^g].$
\item $H^i(\kk[V],\kk[V^g]) \cong \Lambda^i V \otimes \kk[V^g].$
\item The cohomological map $H^i(\kk[V],\kk[V]g) \to H^i(\kk[V],\kk[V^g])$
which arises from the morphism $\kk[V]g \to \kk[V^g]g \cong \kk[V^g]$
is the natural inclusion 
$\Lambda^{i - d_g} V^g \otimes  \Lambda^{d_g}(V^g)^\vee \otimes \kk[V^g] \to \Lambda^i V \otimes \kk[V^g]$.
\end{enumerate}
\end{prop}

{\it Proof.} Choose a basis $\{ v_1, \ldots, v_n\}$ in $V$ so that $v_1, \ldots, v_{n-d_g}$ are
$g$-invariant and span $V^g$, and $v_{n-d_g+1},\ldots, v_n$ are eigenvectors of $g$ with eigenvalues $\lambda_i \ne 1$,
hence span $(V^g)^\vee$. Let $x_i$ be the dual basis of $V^* \subset A$.

The algebra $\kk[V]$ has a Koszul resolution, so the groups
$H^i(\kk[V],\kk[V]g)$ can be computed as the cohomology of the complex
$$
0 \leftarrow \Lambda^n V \otimes \kk[V]g \leftarrow \Lambda^{n-1} V \otimes \kk[V]g \leftarrow \ldots%
\leftarrow V \otimes \kk[V]g \leftarrow \kk[V]g \leftarrow 0
$$
with the differential
$$
d(\xi \otimes ag) =%
\sum\limits_{i=1}^n v_i \wedge \xi \otimes (x_i-x_i^g) ag =%
\sum\limits_{i=n-d_g+1}^n v_i \wedge \xi \otimes (1-\lambda_i)x_iag
$$
and the groups $H^i(\kk[V],\kk[V^g])$ as the cohomology of
$$
0 \leftarrow \Lambda^n V \otimes \kk[V^g] \leftarrow \Lambda^{n-1} V \otimes \kk[V^g] \leftarrow \ldots%
\leftarrow V \otimes \kk[V^g] \leftarrow \kk[V^g] \leftarrow 0
$$
with zero differential.
Then the proposition is straightforward.
$\square$.

The situation of a vector space with a linear action of a finite group has a nice property:
the Koszul complex $\Lambda^*V \otimes \kk[V]g$ is a direct summand in the Hochschild complex
$Hom(\kk[V]^{\otimes *}, \kk[V]g)$, and the complex formed by the cohomology (with
zero differential) $\Lambda^{i - d_g} V^g \otimes  \Lambda^{d_g}(V^g)^\vee \otimes \kk[V^g]$
is in turn a direct summand in the Koszul complex (here we use the fact that $\kk[V^g]$ 
is a submodule of $\kk[V]g$), hence $H^i(\kk[V], \kk[V]g)$ is a direct summand in
$C^i(\kk[V], \kk[V]g)$. Now the map $\mu$ from Lemma \ref{mu}, combined with the projection,
gives a map $H^i(\kk[V], \kk[V]g) \otimes H^i(\kk[V], \kk[V]h) \to H^{i+j}(\kk[V], \kk[V]gh)$
that coincides with the
multiplication in the $G$-invariant part.
In the setting of Proposition \ref{loc_isom} this map may be computed directly.
To do this in a convenient way, we need following lemma:
\begin{lemma}
If $(V^g)^{\vee} \cap (V^h)^{\vee} = \{0\}$, then $V^{gh} = V^g \cap V^h$.
\end{lemma}\label{intersect}
{\it Proof.} Obvoiusly $V^g \cap V^h \subset  V^{gh}$.
The symmetrization by $g$ (denoted by $\pi^g$) projects $V$ onto $V^g$.
Let $(V^g)^{\vee} \cap (V^h)^{\vee} = \{0\}$. The kernel of $\pi^g$ is $(V^g)^{\vee}$. 
If we take any $v\in V^{gh}$, then $v = gh\cdot v$ and
$0=\pi^gv-\pi^ggh\cdot v = \pi^g(v- h\cdot v)$, hence $(Id-h)v \in (V^g)^{\vee}$.
But the operator $Id-h$ is the projection onto $(V^h)^{\vee}$, hence
$(Id-h)v=0$ and $v\in V^h$. By the same argument $v\in V^g$, and since $v$
was an arbitrary element of $V^{gh}$, we have $V^{gh} \subset V^g \cap V^h$.
$\square$

\begin{prop}\label{mult_local}
 The map 
$$
\tilde\mu : \Lambda^{i - d_g} V^g \otimes  \Lambda^{d_g}(V^g)^\vee \otimes%
\Lambda^{j - d_h} V^h \otimes  \Lambda^{d_h}(V^h)^\vee \to%
\Lambda^{i+j - d_{gh}} V^{gh} \otimes  \Lambda^{d_{gh}}(V^{gh})^\vee
$$
which is zero if $V^{gh} \ne V^g \cap V^h$ and sends
$(\xi^1\otimes \xi^2\otimes f) \otimes (\nu^1\otimes \nu^2 \otimes e)$ to
$(-1)^{d_g(j-d_h)} \pi^{gh}\xi^1 \wedge \pi^{gh}\nu^1 \otimes \xi^2\wedge\nu^2 \otimes (fe)|_{V^{gh}}$
otherwise, induces the multiplication on the $G$-invariant part, which
coincides with the cup product.
\end{prop}
{\it Proof.} The maps 
$$\psi_{i,g}: Hom (\kk[V]^{\otimes i}, \kk[V]g) \to  \Lambda^{i - d_g} V^g \otimes  \Lambda^{d_g}(V^g)^\vee \otimes \kk[V^g]$$
may be written down as follows:
\begin{multline*}
\psi_{i,g}(\alpha g) =\\
= \sum\limits_{1\leq m_1 <\ldots \leq n-d_g}
\alpha (x_{m_1} \otimes\ldots x_{m_{i-d_g}} \otimes x_{n-d_g+1}\otimes\ldots x_n)%
v_{m_1}\wedge\ldots v_{m_{i-d_g}} \otimes v_{n-d_g+1}\wedge\ldots v_n f|_{V^g}.
\end{multline*}
As for the maps in the inverse direction, they are quite complicated, but since $\psi_{i,g}$
are computed only on linear functions, we can only consider the linear parts of cochain maps:
$$\phi_{i,g}: \Lambda^i V \otimes \kk[V^g] \to Hom ((V^*)^{\otimes i}, \kk[V]g)$$
which are of the form
$$
\phi_{i,g}(\xi_1 \wedge\ldots\wedge \xi_i \otimes f) (y_1\otimes\ldots\otimes y_i) =
\frac{1}{i!}\sum\limits_{\sigma \in S_i} (-1)^{\sigma}
\langle \xi_1, y_{\sigma(1)}\rangle \cdot\ldots\cdot  \langle \xi_i, y_{\sigma(i)}\rangle fg.
$$

Hence the product $\psi_{i+j,gh} \mu (\phi_{i,g}\otimes\phi_{j,h})$ of two vector fields
$\xi = \xi_1\wedge\ldots\wedge \xi_i \otimes f \in \Lambda^i V \otimes \kk[V^g]$ and
$\nu = \nu_1\wedge\ldots\wedge \nu_j \otimes e \in \Lambda^i V \otimes \kk[V^h]$ is
the $\Lambda^{i+j - d_{gh}} V^{gh} \otimes  \Lambda^{d_{gh}}(V^{gh})^\vee$-component
of $\xi_1\wedge\ldots\wedge \xi_i \wedge \nu_1^g \wedge\ldots\wedge \nu_j^g \otimes (fe^g)|_{V^{gh}}
\in \Lambda^i V \otimes \kk[V^{gh}]$. 
Then if $\xi$ and $\nu$ are taken in the cohomology,
$\xi$ contains $\Lambda^{top} (V^g)^{\vee}$, and $\nu$ contains $\Lambda^{top} (V^h)^{\vee}$.
If $(V^g)^{\vee} \cap (V^h)^{\vee} \ne \{0\}$, then $(V^g)^{\vee} \cap ((V^h)^{\vee})^g \ne \{0\}$,
and the product is automatically zero; if $(V^g)^{\vee} \cap (V^h)^{\vee} = \{0\}$, then by Lemma \ref{intersect}
$V^{gh} = V^g \cap V^h$ and $(V^{gh})^{\vee} = (V^g)^{\vee} \oplus (V^h)^{\vee} = (V^g)^{\vee} \oplus ((V^h)^{\vee})^g$.
Note that in this case $(fe^g)|_{V^{gh}} = (fe)|_{V^{gh}}$.
Now its time to write $\xi = \xi^1 \otimes \xi^2 \in \Lambda^{i - d_g} V^g \otimes  \Lambda^{d_g}(V^g)^\vee$,
$\nu = \nu^1 \otimes \nu^2 \in \Lambda^{i - d_h} V^h \otimes  \Lambda^{d_h}(V^h)^\vee$.
Then the product may be rewritten as $(-1)^{d_g(j-d_h)}\xi^1 \wedge (\nu^1)^g \wedge \xi^2 \wedge (\nu^2)^g \otimes (fe)|_{V^{gh}}$.
Note that for any vector $w \in V$ we have $w-w^g \in (V^g)^{\vee}$, so after taking 
the wedge product with $\xi^2 \in \Lambda^{top}(V^g)^{\vee}$ we have
$w\wedge\xi^2= w^g\wedge \xi^2$, hence the product simply equals
$(-1)^{d_g(j-d_h)}\xi^1 \wedge \nu^1 \wedge \xi^2 \wedge \nu^2 \otimes (fe)|_{V^{gh}}$, which
in turn equals $\pi^{gh}\xi^1 \wedge \pi^{gh}\nu^1 \otimes \xi^2\wedge\nu^2 \otimes fe|_{V^{gh}}$,
and we are done.
$\square$.

%
%

If $V$ carries a symplectic form $\omega$, and the action of $G$ is
symplectic, then for all $g$ the forms $\omega|_{(V^g)^\vee}$ 
are nondegenerate, and we can construct nonzero elements
$s_g \in \Lambda^{d_g} (V^g)^\vee$ such that if 
$V^{gh} = V^g \cap V^h$, then $s_{gh} = s_g \wedge s_h$;
namely, take $s_g$ dual to $(\omega|_{(V^g)^\vee})^{\wedge d_g/2}$.
Using sections $s_g$, another
canonical isomorphism can be established:
\begin{prop}\label{mult_sympl} If $V$ is symplectic, then
$H^i(\kk[V],\kk[V]g) \cong \Lambda^{i - d_g} V^g \otimes \kk[V^g]$. The multiplication
is given by the wedge product: the product of two vector fields
$\xi_g \in \Lambda^{i - d_g} V^g \otimes \kk[V^g]$ and $\xi_h \in \Lambda^{i - d_h} V^h \otimes \kk[V^h]$
is zero if $V^{gh} \ne V^g \cap V^h$ and $(-1)^{d_g (j-d_h)} \xi_g \wedge \xi_h$ otherwise.
\end{prop}

\section{The affine case}

Now we can return to the case of smooth affine algebraic $X$
and $A=\kk[X]$ its algebra of regular functions. The group $G$ acts
algebraically on $X$; for $g\in G$ denote by $X^g_m$ the 
$m$-th connected component of the
algebraic subvariety of $g$-invariant points. By finiteness of $G$
the variety $X^g_m$ is smooth. Denote by $A_{g, m}=\kk[X^g_m]$ its algebra
of regular functions.

\begin{lemma}\label{simeq} $A_{g, m} \simeq A_{g, m}g$ as $A$-bimodules.
\end{lemma}

Note that the normal bundle $N_XX^g_m$ is naturally embedded into the
restriction $TX|_{X^g_m}$ as a subbundle generated by $g$-semiinvariant vectors.

\begin{prop}\label{global_isom}
\begin{enumerate}
\item $H^i(A, Ag) \simeq \bigoplus\limits_m \Gamma (X^g_m, \Lambda^{i-d_{g, m}} TX^g_m \otimes \Lambda^{d_{g, m}} N_XX^g_m)$;
\item $H^i(A, A_{g, m}) \simeq \Gamma (X^g_m, \Lambda^i TX|_{X^g_m})$;
\item The cohomological map $H^i(A, Ag) \to \bigoplus\limits_m H^i(A, A_{g, m})$ induced by the map
$Ag \to \bigoplus\limits_m A_{g, m}g \simeq \bigoplus\limits_m A_{g, m}$ comes from the natural inclusions
$\Lambda^{i-d_{g, m}} TX^g_m \otimes \Lambda^{d_{g, m}} N_XX^g_m \to \Lambda^i TX|_{X^g_m}$.
\end{enumerate}
\end{prop}

{\it Proof.} It is convenient to switch to the language of coherent sheaves.
By definition $H^i(A, \cdot) = Ext^i_{A \otimes A^{op}}(A, \cdot)$; the algebra $A$
is commutative, hence $A^{op}\cong A$, and $A \otimes A = \Gamma (X\times X, \Oo_{X\times X})$.
Let $\Delta \subset X\times X$ be the diagonal, $Y \subset X\times X$ be the
graph of $g: X \to X$.Then $H^i(A, Ag) = Ext^i_{X \times X} (\Oo_{\Delta}, \Oo_Y)$,
and $H^i(A, A_{g, m}) = Ext^i_{X \times X} (\Oo_{\Delta}, \Oo_{\Delta^g_m})$.

Our main tool in the proof of the proposition will be the following lemma:
\begin{lemma}\label{formal}
 Let $F, G$ be coherent sheaves on an algebraic variety $Y$,
let $f: F \to G$ a morphism of coherent sheaves.
Then $f$ is an isomorphism (resp. injective) iff for any point $y \in Y$ in a formal
neighborhood of $y$ it becomes an isomorphism (resp. injective).
\end{lemma}

Then the proof goes as follows: we construct the maps of sheaves on $X\times X$:
\begin{multline}\label{phi}
Ext^i_{X \times X} (\Oo_{\Delta}, \Oo_Y)
\xrightarrow{ \phi'}
\bigoplus\limits_m Ext^i_{X \times X} (\Oo_{\Delta}, \Oo_{\Delta^g_m})\\
\bigoplus\limits_m Ext^i_{X \times X} (\Oo_{\Delta}, \Oo_{\Delta^g_m})
\xleftarrow{\phi}
\bigoplus\limits_m \Lambda^i T\Delta|_{\Delta^g_m} \\
\bigoplus\limits_m \Lambda^i T\Delta|_{\Delta^g_m}
\xrightarrow{\phi''}
\bigoplus\limits_m \Lambda^{i-d_{g, m}} T\Delta^g_m \otimes \Lambda^{d_{g, m}} N_{\Delta}\Delta^g_m
\end{multline}
then apply Lemma \ref{formal} to prove that $\phi$ is an isomorphism and
$\phi'$ is an inclusion; then we can take the composition map
$\psi = \phi'' \circ \phi^{-1} \circ \phi'$ and after applying Lemma \ref{formal} once more,
see that $\psi$ is an isomorphism.

Let us define the maps: $\phi'$ is the functorial map arising from
$\Oo_Y \to \oplus \Oo_{\Delta^g_m}$ (note that $\cup \Delta^g_m = \Delta\cap Y$);
the last one, $\phi''$, is the projection (recall that
$T\Delta|_{\Delta^g_m} = T\Delta^g_m \oplus N_{\Delta}\Delta^g_m$).
To define $\phi$, take a resolution of $\Oo_{\Delta}$ by free $\Oo_{X \times X}$-sheaves
$\Gamma(X \times X, \Oo_{\Delta})^{\otimes n} \otimes \Oo_{X \times X}$
(in fact, the bar resolution) and construct the map 
$\phi_i: \Lambda^i T\Delta |_{\Delta^g_m} \to
Hom_{\kk} (\Gamma(X \times X, \Oo_{\Delta})^{\otimes i}, \Gamma(X \times X, \Oo_{\Delta^g_m}))$
directly:
\begin{equation}\label{explicit}
(\phi_i (\xi_1 \wedge\ldots\wedge \xi_i))(a_1 \otimes\ldots\otimes a_i) =
\frac{1}{i!}
\sum\limits_{\sigma \in S_i} (-1)^{\sigma} \partial_{\xi_1} a_{\sigma(1)} \cdot\ldots\cdot  \partial_{\xi_i} a_{\sigma(i)},
\end{equation}
where RHS is a well-defined function on $\Delta^g_m$.
By an easy calculation, the image of $\phi$ lies in the kernel of the differential
in the bar resolution, so it induces a map
$\phi: \Lambda^i T\Delta |_{\Delta^g_m} \to Ext_{X\times X} (\Oo_{\Delta}, \Oo_{\Delta^g_m})$.

When $g=1$ this map is the famous Hochschild-Kostant-Rosenberg isomorphism
\cite{HKR}, which appears in this exact form in \cite{Kontsevich}, 4.6.1.1.

To proceed, we need one more lemma:
\begin{lemma}\label{polydisc}
If the order $|G|$ of a finite group $G$ is prime to the
characteristic of the ground field $\kk$, then any action of $G$ on a formal polydisc over $\kk$
is equivalent to a linear action.
\end{lemma}

Now for $x \in X \times X - \Delta^g$ all the considered sheaves are zero in the formal
neighborhood of $x$, and for $x \in \Delta^g_m$, which corresponds to a point
$x \in X^g_m$, the action of $G$ on the formal neighborhood
of $x \in X$ is linear by Lemma \ref{polydisc}, and we are in the situation of Proposition \ref{loc_isom} with $V = T_xX$.
All previous constructions of cohomology groups and their maps commute
with the transition to the formal neighborhood, and in the formal neighborhood
the statements follow from Proposition \ref{loc_isom}. Then we can apply Lemma \ref{formal}
to finish the proof.
$\square$

Propositions \ref{A_Ag} and \ref{global_isom} together prove Theorem \ref{isom}.

Note that (\ref{explicit}) actually defines a map of $\Gamma(X^g_m, \Lambda^i TX |_{X^g_m})$
into Hochschild cochains $C^i(A, A_{g, m})$, so we can introduce a multiplication:
$H^i(A,Ag)\otimes H^j(A,Ah) \to C^{i+j}(A, Agh)$ 
The equalities from
Theorems \ref{mult} and \ref{main_sympl} hold locally by Propositions \ref{mult_local} and
\ref{mult_sympl} hence they hold globally.

\vspace{2cm}
\noindent Department of Mathematics, Harvard University, 1 Oxford Street, Cambridge, MA 02138, USA.

\noindent email: {\tt rina@math.harvard.edu}


\begin{thebibliography}{}

\bibitem{Caldararu}
A. Caldararu, {\it The Mukai pairing, II: the Hochschild-Kostant-Rosenberg isomorphism,}
math.AG/0308080.

\bibitem{DE}
V. Dolgushev, P. Etingof, {\it Hochschild cohomology of quantized 
symplectic orbifolds and the Chen-Ruan cohomology,}
Int. Math. Res. Not. 2005, no. {\bf 27}, 1657-1688, also math.QA/0410562.

\bibitem{GK}
V. Ginzburg, D. Kaledin, {\it Poisson deformations of symplectic
quotient singularities,} Adv. Math. {\bf 186}, 1 (2004) 1-57,  
also math.AG/0212279.

\bibitem{HKR}
G. Hochschild, B. Kostant, A. Rosenberg, {\it Differential forms on
regular affine algebras,} Trans. AMS {\bf 102} (1962), 383-408.

\bibitem{Kaledin}
D. Kaledin, {\it Multiplicative McKay correspondence in the symplectic case,}
math.AG/0311409.

\bibitem{Kontsevich}
M. Kontsevich, {\it Deformation quantization of Poisson manifolds, I,}
Lett. Math. Phys. {\bf 66} 157-216 (2003), also q-alg/9709040.

\bibitem{NPPT}
N. Neumaier, M. J. Pflaum, H. B. Posthuma, X. Tang,
{\it Homology of formal deformations of proper etale Lie groupoids},
math.KT/0412462.

\bibitem{Witherspoon}
S. J. Witherspoon, {\it Products in Hochschild cohomology and Grothendieck rings
of group crossed products,} math.RA/0212003.

\end{thebibliography}
\end{document}